\documentclass[11pt]{amsart}
\usepackage{geometry}                
\usepackage{graphicx}
\usepackage{amssymb}
\usepackage{epstopdf}
\usepackage{hyperref}
\newtheorem{thm}{Theorem}

\title[Tossing half-coins and other partial coins]{Tossing half-coins and other partial coins: \\ signed probabilities and Sibuya distribution}

\date{}                                

\author{Nikolai Leonenko}
\address{Nikolai Leonenko\\ 
		School of Mathematics \\
		Cardiff University, UK}
\email{leonenkon@cardiff.ac.uk} 

\author{Igor Podlubny}
\address{Igor Podlubny\\ 
		BERG Faculty \\ 
		Technical University of Kosice, Slovakia}
\email{igor.podlubny@tuke.sk}

\begin{document}

\markboth{\uppercase{Nikolai Leonenko, Igor Podlubny}}
{ \uppercase{Tossing half-coins and other partial coins \ldots}}

\thispagestyle{empty}

\begin{abstract}
	A method for numerical simulation of signed probability
	 distributions for the case of tossing $1/n$-th of a coin is presented
	 and illustrated by examples. 
\end{abstract}

\maketitle

\section{Introduction}

What do fractional-order differentiation~\cite{LP-MCFD-1,LP-MCFD-2,LP-MCFD-3} on one hand, and ``partial coins''~\cite{Szekely-2005}, on the other, have in common? 
The answer is: the need to simulate signed probability distributions, and the use of Sibuya probability distribution. 

Recently, we succeeded in developing the Monte Carlo method 
for numerical fractional-order differentiation~\cite{LP-MCFD-1}, 
and the key tool was the Sibuya distribution~\cite{Sibuya-1979}.
For orders of differentiation higher than one, this led us to working 
with signed probability distributions, where the number of negative probabilities was finite 
and the number of positive probabilities was infinite (or vice versa), 
and we showed how to deal with such situations~\cite{LP-MCFD-2, LP-MCFD-3}. 

In the present paper we  provide -- to our best knowledge -- the first examples
of numerical simulation of signed probability distributions related to random variables 
that can be called \textit{partial coins}, and in this case the numbers of negative and positive 
probabilities are both infinite, which did not allow the use of our previous method. 
However, the probability distributions of \textit{partial coins} appeared to be 
related to the Sibuya  distribution.

Below we start with explaining what are half-coins and other partial coins, 
and how they are related to signed probability distributions. 
Then we turn to the Sibuya probability distribution and to our methods of its simulation. 
After that, we use the Sibuya probability distribution for obtaining 
the decomposition of the partial coin probability distributions 
into two probability distributions with only positive probabilities. 
Finally, we provide examples of numerical simulation of partial coins. 

\section{Half-coins and other partial coins}

A classical fair coin is a random variable $X$ that takes the values 0 and 1 with probability 1/2, and the expectation equal to 1/2. Its probability generating function is $c(x) =(1+x)/2$. The addition of independent random variables corresponds to multiplication of their probability generating functions. Therefore, the probability generating function of the sum of two fair coins is equal to $((1+x)/2)^2$.

\textit{Half of a coin}, or more conveniently \textit{a half-coin}, was introduced by Sz\'ekely \cite{Szekely-2005}, who noticed that
it is natural to define a half-coin via its probability generating function
\begin{equation}\label{eq:half-coin}
	\left(\frac{1+x}{2}\right)^\mu = \sum_{n=0}^{\infty} p_n x^n, 
	\quad 
	\mu = \frac{1}{2}.
\end{equation}
where the coefficients $p_n = 2^{-\mu} {\mu \choose n}$ have alternating signs starting from $n=1$:
\begin{equation} \label{eq:alternating-signs}
	p_0 > 0, \quad
	p_1 > 0, \quad
	p_2 < 0, \quad
	p_3 > 0, \quad
	p_4 < 0, \quad
	p_5 > 0, \quad
	p_6 < 0, \quad \ldots
\end{equation}

A half-coin corresponds to \hbox{$\mu=1/2$}, and it 
takes the values $n$ with signed probabilities~-- positive for $n=0, 1, 3, 5, \ldots$,  and negative for $n = 2, 4, 6, \ldots$. 
The repeated ``flipping'' of two half-coins gives the expectation equal to 0.5, 
similarly to the case of the repeated flipping of one normal coin.  

Taking $\mu = 1/3$ yields one-third-coins and the repeated ``flipping'' of three one-third-coins makes a normal coin with  the expectation equal to 0.5; $\mu = 1/4$ defines quarter-coins,  the repeated flipping of four quarter-coins makes a normal coin, etc. 

Non-integer values of $\mu$ produce $\mu$th-coins, and to get a normal coin it is necessary to ``flip''  $1/\mu$ of such $\mu$th-coins, 
which depending on $\mu$ can be a finite or even infinite number of such  \textit{partial coins}.

Obviously, from (\ref{eq:half-coin}) follows that 
$$
	\sum_{n=0}^{\infty} p_n = 1, 
	\qquad 
	\sum_{n=0}^{\infty} |p_n| = 2^{1 - \mu} .
$$

For the case of  $\mu = 1/2$, Sz\'ekely obtained an expression for the coefficients $p_n$ in terms of the Catalan numbers $C_n$:
$$
	p_n=(-1)^{n-1} \sqrt{2} \, \frac{C_{n-1}}{4^n}, \quad n = 0, 1, 2, \ldots
$$ 
where
$$
	C_n = \frac{{2n \choose n}}{n+1}, \quad n = 0, 1, 2, \ldots; 
	\qquad
	C_{-1} = -\frac{1}{2}. 
$$

\section{Sibuya probability distribution}

The probability distribution introduced by Sibuya~\cite{Sibuya-1979}
and studied first by Sibuya and Shimitzu~\cite{Sibuya-Shimitzu-1981} 
and then, with some generalizations, by other authors \cite{Devroye-1993,Pillai-Jayakumar-1995,KP-2018}, 
can be introduced in the following way.

Let us consider the sequence of independent Bernoulli trials, in which the $k$-th trial 
has probability of success $\alpha/k$ ($0<\alpha<1$, \, $k = 1, 2, \ldots$), 
and let $Y \in \{1, 2, \ldots \}$ be the trial number in which the first success occurs.
Then 
$$
	\mathbb{P} (Y=1) = p_1 = p_1(\alpha) = \alpha,
$$
$$
	 \mathbb{P} (Y=k) = p_k = p_k(\alpha) =
	 (1-\alpha) 
	 (1-\frac{\alpha}{2})
	\ldots
	(1-\frac{\alpha}{k-1}) 
	\frac{\alpha}{k},
	\quad k \geq 2,
$$	
or, in other words, its probability mass function is
\begin{eqnarray*}
	p_k = \mathbb{P} (Y=k) & = &
	(-1)^{k+1}
	\frac{\alpha (\alpha -1) \ldots (\alpha - k + 1)}{k!} = \nonumber\\
	& = &
	 (-1)^{k+1} { \alpha \choose k} , \qquad k = 1,2, \ldots , \label{S-mass-function}
\end{eqnarray*}
and its cumulative distribution function is
\begin{eqnarray*}
	F_k = \sum_{j=1}^{k} p_j 
	& = & 
	1 - (-1)^{k} {\alpha -1  \choose k } = 
	1- \frac{1}{k \, B (k, 1-\alpha)}  =  \nonumber\\
	& = & 
	1 - \frac{\Gamma (k - \alpha + 1)}{k \, \Gamma (k) \Gamma (1-\alpha)},
	\quad
	k = 1,2, \ldots , \label{S-cumulative-function}
 \end{eqnarray*}
where $B(x,y)$ and  $\Gamma(x)$  are Euler's beta and gamma functions.

The probability generating function of the Sibuya distribution is 
\begin{eqnarray} 
	G_Y (x) & = & \mathbb{E} \, x^{Y} =
 	\sum_{k=1}^{\infty}
	x^{k}
	(-1)^{k+1} {\alpha \choose k}  \\
	& = & 
	\sum_{k=1}^{\infty}
	x^{k} (-1)^{k+1} w_k^{(\alpha)}  
	= 
	1 - (1-x)^\alpha, 
	\quad 
	|x| < 1.   \label{eq:Sibuya-GF}
\end{eqnarray}

\noindent
where
\begin{equation}
	w_k^{(\alpha)} = { \alpha \choose k}.
\end{equation}

For $0 < \alpha < 1$, all coefficients 
in the power series expansion 
of the generating function (\ref{eq:Sibuya-GF}) are positive, monotonically decreasing, 
and their sum is equal to one.

\section{Simulation of the signed distribution}

It is known that summation of independent random variables corresponds to the product of their generating functions. 
The way to the simulation of signed probability distributions was indicated by Székely~\cite{Szekely-2005}:

\begin{thm}[G. Sz\'ekely~\cite{Szekely-2005}]\label{theorem-1}
For every generalized generating function $f$ of a signed probability distribution 
there exist two generating functions $g$ and $h$ of ordinary non-negative probability distributions
such that $fg =h$.
\end{thm}

The proof of this theorem can be found in \cite[Theorem 1]{Ruzsa-Szekely-1983} 
or in \cite[Lemma 3.6]{Ruzsa-Szekely-APT}.

This theorem guarantees the existence of $g$ and $h$, but not the uniqueness. 
We need a pair of distributions $g$ and $h$, the difference of which gives the signed distribution $f$.\footnote{This is similar to the situation 
with vectors: $a = b - c = b_1 -c_1$, where
\begin{itemize}
 	\item[] $a = [2,  -1,  3, 5]$,   with positive and negative components, 
	\item[] $b = [5,  6,  7, 8]$,   $c = [3, 7, 4, 3]$, both vectors with only positive components, 
	\item[] $b_1 = [7,  8,  8, 10]$,   $c_1= [5, 9, 5, 5]$, both vectors with only positive components.
\end{itemize}}

The Theorem~\ref{theorem-1} does not provide any method 
for constructing the functions $g$ and $h$. However, in the case of partial coins
we, based on our experience~\cite{LP-MCFD-1,LP-MCFD-2,LP-MCFD-3},
succeeded in guessing these functions.   

For the $1/n$-coin the generating function of its signed probability distribution $f$ is:
\begin{equation} \label{eq:f}
	f (x) = \left(\frac{1+x}{2}\right)^{1/n}
\end{equation}

The pairs of suitable ordinary nonnegative probability distributions $g$ and $h$ are:
\begin{equation} \label{eq:g}
	g_k(x) = \Bigl( 1 - (1-x)^{1/n} \Bigr) x^k,
\end{equation}
\begin{equation} \label{eq:h}
	h_k(x) =  2^{1/n} x^k \Bigl( (1+x)^{1/n} - (1-x^2)^{1/n}  \Bigr) ,
\end{equation}
where obviously $f(x) g_k(x) = h_k (x)$, and $k$ can be $-1, 0, 1, 2, 3, 4, \ldots$.

The functions $g_k(x)$ are given by the generating function of the Sibuya distribution multiplied by $x^k$,
and $h_k(x)$ are the products of $f(x)$ and $g_k(x)$. 

The pairs of distributions $g_k$ and $h_k$ with the support on $\{k+1, k+ 2, \ldots  \}$ work. 
%
\begin{equation}
	f g_k = h_k
\end{equation}
\begin{equation} \label{eq:Hi-Gi}
	F_i = H_i - G_i
\end{equation}

If $k = -1$, then $g$ and $h$  have the same support $\{ 0, 1, 2, 3, \ldots \}$

\section{Coefficients of the generating functions $f$, $g$, and $h$}

Let us recall that the coefficients of the power series expansion of a probability generating function 
mean probabilities with which the random variable, say $Y$, take the values from its support. 

Let us denote $\alpha = 1/n$. 
For $0 < \alpha < 1$ the expansion of the generating function 
\begin{equation}
	f (x) = \left(\frac{1+x}{2}\right)^{\alpha}
\end{equation}
into a power series produces coefficients with alternating signs. 

The coefficients of the power series expansion of $g(x)$ for \hbox{$0<\alpha<1$} are all positive:
\begin{equation}
	g_k(x) = \Bigl( 1 - (1-x)^{\alpha} \Bigr) x^k.
\end{equation}

We have to verify that the coefficients of the power series expansion of 
\begin{equation} \label{eq:hk}
	h_k(x) = f(x) g_k(x) =  2^{-\alpha} x^k \Bigl( (1+x)^{\alpha} - (1-x^2)^{\alpha}  \Bigr)
\end{equation}
are also all positive. 

For this, we have to examine the power series expansion of the expression
\begin{equation} \label{eq:D}
	D(x) = (1+x)^{\alpha} - (1-x^2)^{\alpha}
\end{equation}

We have:
\begin{eqnarray}
	D(x)  & = & 
		\sum_{m=1}^{\infty}  { \alpha \choose m} x^m 
		- \sum_{n=1}^{\infty} (-1)^n { \alpha \choose n} x^{2n} \\
	& = & 
		\sum_{r = 0}^{\infty} { \alpha \choose 2r+1} x^{2r+1} 
		+ \sum_{r = 1}^{\infty} { \alpha \choose 2r} x^{2r}
		-  \sum_{n=1}^{\infty} (-1)^n { \alpha \choose n} x^{2n} \\
	& = & 
		\sum_{r = 0}^{\infty} { \alpha \choose 2r+1} x^{2r+1} 
		+ 
		\sum_{n=1}^{\infty} \Bigl\{ { \alpha \choose 2n} - (-1)^n {\alpha \choose n}   \Bigr\} x^{2n}  \\
	& = & 
		\sum_{r = 0}^{\infty} w_{2r+1}^{(\alpha)}  x^{2r+1} 
		+ 
		\sum_{n=1}^{\infty} \Bigl\{ w_{2n}^{(\alpha)} - (-1)^n w_{n}^{(\alpha)}   \Bigr\} x^{2n}  . 
		\label{eq:two-sums}
\end{eqnarray}

For \hbox{$0<\alpha<1$}, the binomial coefficients $w_k^{(\alpha)}$ have alternating signs, 
and the coefficients at the odd powers of $x$ are positive, 
so coefficients of the first sum in (\ref{eq:two-sums}) are all positive.
We need to verify that the coefficients in the second sum are also positive.

Let us denote 
$$
	b_k^{(\alpha)} = | w_k^{(\alpha)} | = \left|  { \alpha \choose k}  \right|,
$$
then 
$$
	w_k^{(\alpha)} =  { \alpha \choose k} = (-1)^{k-1} b_k^{(\alpha)}. 
$$

This allows us to write
\begin{eqnarray*}
	\Bigl\{ w_{2n}^{(\alpha)} - (-1)^n w_{n}^{(\alpha)}   \Bigr\}  
	& = & 
	(-1)^{2n-1} b_{2n}^{(\alpha)}  - (-1)^n \cdot (-1)^{n-1} b_n^{(\alpha)} \\
	& = & 
	(-1)^{2n-1} b_{2n}^{(\alpha)}  + (-1)^{2n}  b_n^{(\alpha)}  \\
	& = & 
	b_n^{(\alpha)} - b_{2n}^{(\alpha)}. 
\end{eqnarray*}

For \hbox{$0<\alpha<1$}, the sequence of positive numbers $b_n^{(\alpha)}$ ($n = 1, 2, \ldots$)
is monotonically decreasing. Indeed, using the recurrence relation for the binomial coefficients 
(see \cite[eq. (7.23)]{Podlubny-FDE}), we have:
\begin{eqnarray}
	\frac{b_{n+1}^{(\alpha)}}{b_n^{(\alpha)}}
	& = & 
		\left|
			\frac{w_{n+1}^{(\alpha)}}{w_n^{(\alpha)}}
		\right| 
	= 
		1 - \frac{\alpha + 1}{k+1}  < 1
	   ,
	   \quad
	   k = 1, 2, 3. \ldots
\end{eqnarray}

Therefore, $b_n^{(\alpha)} > b_{n+1}^{(\alpha)}$, and  $b_n^{(\alpha)} > b_{2n}^{(\alpha)}$, 
which means that all coefficients in (\ref{eq:two-sums}) are positive, 
and all coefficients in the power series expansion of $h_k(x)$ are positive.
Also, it follows from (\ref{eq:h}) that the sum of those coefficients is equal to one. 

The above means that we have a family pairs of non-negative probability distributions 
with the generating functions $g_k(x)$ and $h_k(x)$ given by (\ref{eq:g}) and (\ref{eq:h})
which 
can be used for simulating the signed probability distribution 
with the generating function $f(x)$ given by (\ref{eq:f}).

\section{Numerical simulations}

\subsection{Evaluation of the coefficients of $f(x)$ and $g(x)$}
\hspace*{\fill} \\

\noindent
First, let us consider the generating function of the Sibuya distribution (\ref{eq:Sibuya-GF}).
The coefficients of the power series expansion of $g(x)$ for \hbox{$0<\alpha<1$} are all positive:
\begin{equation}
	g_k(x) = \Bigl( 1 - (1-x)^{\alpha} \Bigr) x^k = 
	\sum_{n=1}^{\infty}
	\tilde{w}_n^{(\alpha)}
	x^n, 
	\quad
	\tilde{w}_n^{(\alpha)} >0, 
	\quad 
	n = 1, 2, 3, \ldots \, ,
\end{equation}
where
\begin{equation}
	\tilde{w}_n^{(\alpha)} = (-1)^{n+1} {\alpha \choose n}, 
	\quad 
	n = 1, 2, 3, \ldots
\end{equation}

\begin{figure}
\includegraphics[width=0.9\textwidth]{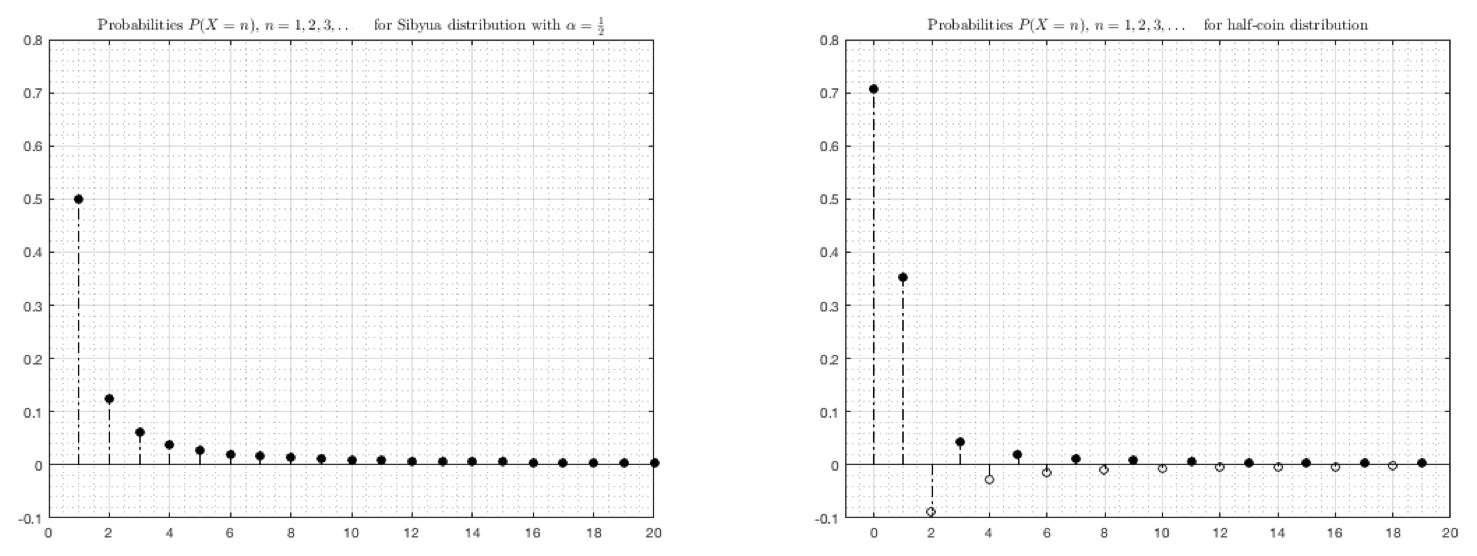}
\caption{Coefficients of $g(x)$ and $f(x)$} 
\label{fig:fg-coefficients}
\end{figure}

The coefficients $\tilde{w}_n^{(\alpha)}$ can be conveniently computed using
 the recurrence relation for the binomial coefficients 
(see \cite[eq. (7.23)]{Podlubny-FDE}). 
For example, for $\alpha=\frac{1}{2}$ the results are shown 
in Figure~\ref{fig:fg-coefficients} on the left. 

The coefficients of the power series expansion 
of the partial-coin distribution (\ref{eq:half-coin}) 
are then computed as
\begin{equation}
	p_0^{(\alpha)} = 2^{-\alpha}, 
	\quad
	p_n^{(\alpha)} = (-1)^{n+1}  \tilde{w}_n^{(\alpha)}, 
	\quad
	n= 1, 2, 3, \ldots
\end{equation}

For $\alpha=\frac{1}{2}$ the results are shown 
in Figure~\ref{fig:fg-coefficients} on the right. 
Starting from $n=2$, the coefficients $p_n$, which are the probabilities $P(X=n)$, 
have alternating signs:
$$
P(X = 2m) < 0, \quad P(X=2m-1)> 0, \qquad m=1, 2, 3 \ldots
$$

\subsection{Evaluation of the coefficients of $h(x)$}
\hspace*{\fill} \\

\noindent
The coefficients $q_k$ of the power series expansion of $h(x)$ 
can be computed using (\ref{eq:h}), (\ref{eq:D}), and (\ref{eq:two-sums}):
\begin{equation}
	q_k^{(\alpha)} = 
	2^{-\alpha}
	\left\{
		\begin{array}{ll}
			p_{2m-1}, & k = 2m-1 \\
			\tilde{w}_m^{(\alpha)} - p_{2m}^{(\alpha)}, & k = 2m \\
		\end{array}
	\right.
\end{equation}

It should be mentioned that expanding $f(x)$, $g(x)$, and $h(x)$ in  power series
using a reliable software for symbolic computations also produces the necessary coefficients. 

\subsection{Method of simulation}\label{sec:simulation}
\hspace*{\fill} \\

\begin{itemize}
	\item[Step 1.]
		Compute the probabilities $\tilde{w}_k^{(\alpha)}$ 
		and $q_k^{(\alpha)}$, $k = 1, 2, 3 \ldots$
	\item[Step 2.]
		Compute the cumulated mass distribution functions $G(x)$	 and $H(x)$.
	\item[Step 3.]
		Generate a set of uniformly distributed points $U_i$ in $(0,1)$.
	\item[Step 4.]
		For each $U_i$  find a pair of the values $G_i$ and $H_i$; 
		for this, the method of the inverse cumulated mass probability 
		function can be used. 
	\item[Step 5.]
		The differences of $H_i$ and $G_i$ give the values $F_i = H_i - G_i$
		of the simulated signed distribution $f$. 
\end{itemize}

\noindent
This approach can be used for simulating the difference of 
any two probability distributions given by their probability mass functions.

\subsection{Examples of numerical simulations of partial coins}
\hspace*{\fill} \\

\noindent
For providing a tool for numerical simulation, a dedicated MATLAB toolbox  
has been created~\cite{Podlubny-partial-coin-toolbox}.  
The function \texttt{partialcoin} has four input parameters 
(type of a coin, number of terms in the expansion of the generating function, 
number of flips, power of the factor), 
and returns 
the vectors of results of flips in the order of their appearance,
the total counts and the frequencies of ones and zeros, 
the expectation of the tossed partial coin, 
and the expectation of such number of the tossed partial coins 
that make a classical fair coin.

\begin{figure}[p!]
\includegraphics[width=0.8\textwidth,height=0.24\textheight]{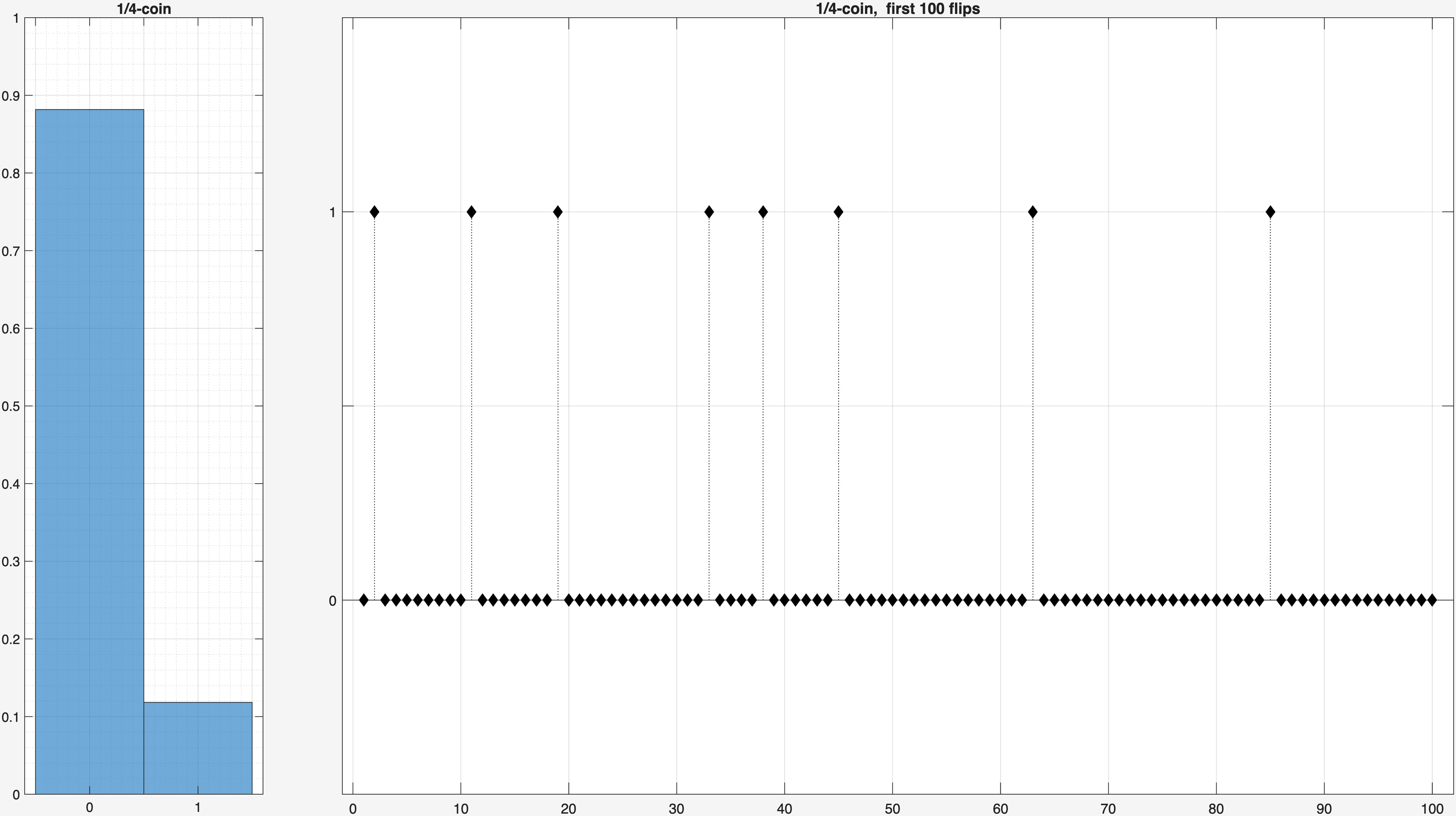}
\caption{A quarter-coin.  \\
Flips: 10000; 
ones: 1183; 
zeros:  8817;  
expectation: 0.1183}
\label{fig:quartercoin}
\end{figure}

\begin{figure}[p!]
\includegraphics[width=0.8\textwidth,height=0.24\textheight]{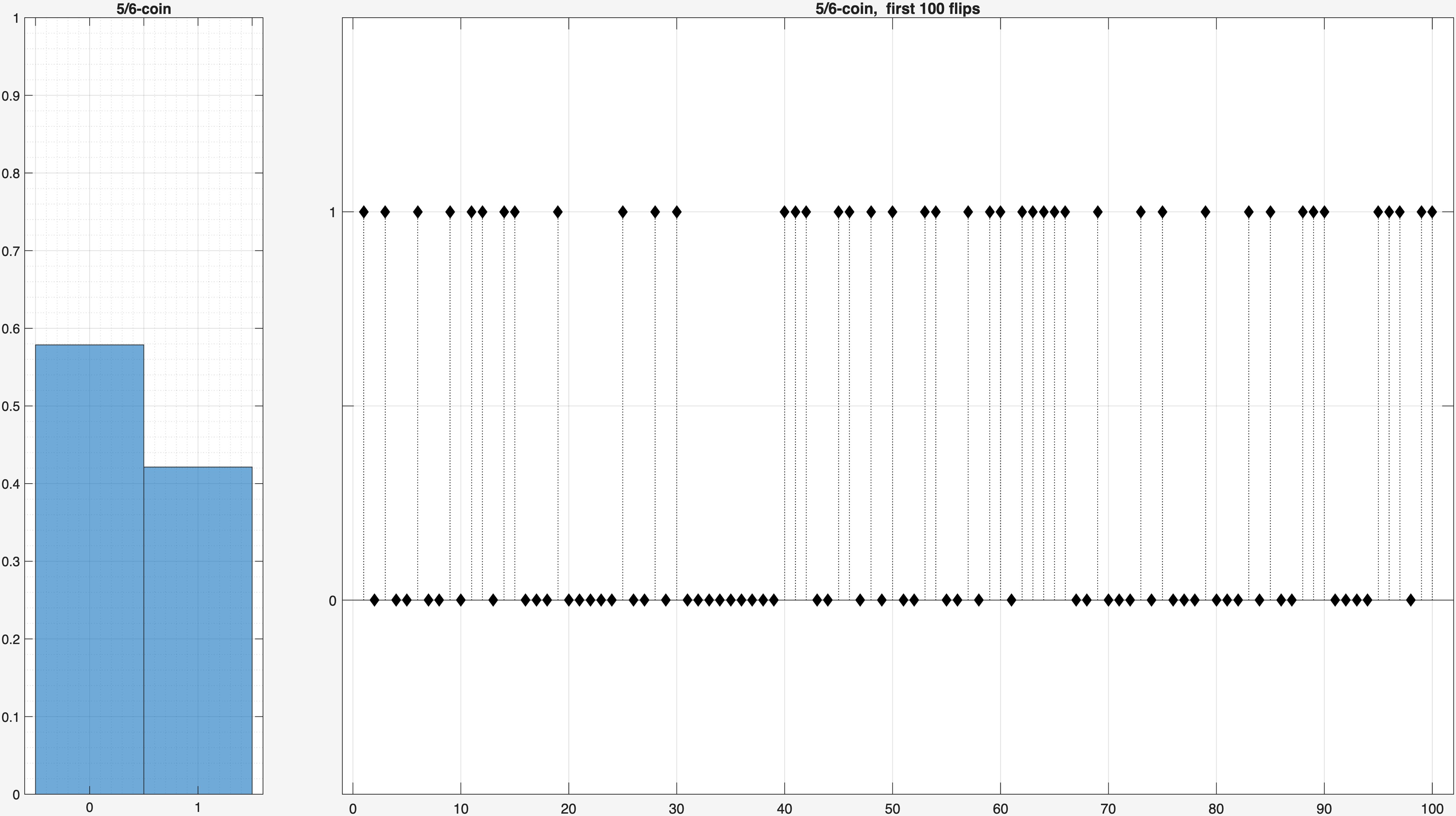}
\caption{A five-sixths-coin. \\
Flips: 10000; 
ones: 4213; 
zeros:  5787;  
expectation: 0.4213} 
\label{fig:five-sixths-coin}
\end{figure}

In Figures~\ref{fig:quartercoin}--\ref{fig:five-sixths-coin} the results 
of numerical simulations are shown for 
a quarter-coin, a one-third-coin, a half-coin, 
a three-quarters-coin, a four-fifths-coin, 
a four-fifths-coin, and a five-sixths-coin. 

The output $F_i$ of each individual flip of a partial coin 
is computed using (\ref{eq:Hi-Gi}), and each flip produces either zero or one.  
The results of the first one hundred flips are shown 
in the right  of Figures~\ref{fig:quartercoin}--\ref{fig:five-sixths-coin},
and the histograms are shown on the left.

\begin{figure}[p!]
\includegraphics[width=0.8\textwidth,height=0.24\textheight]{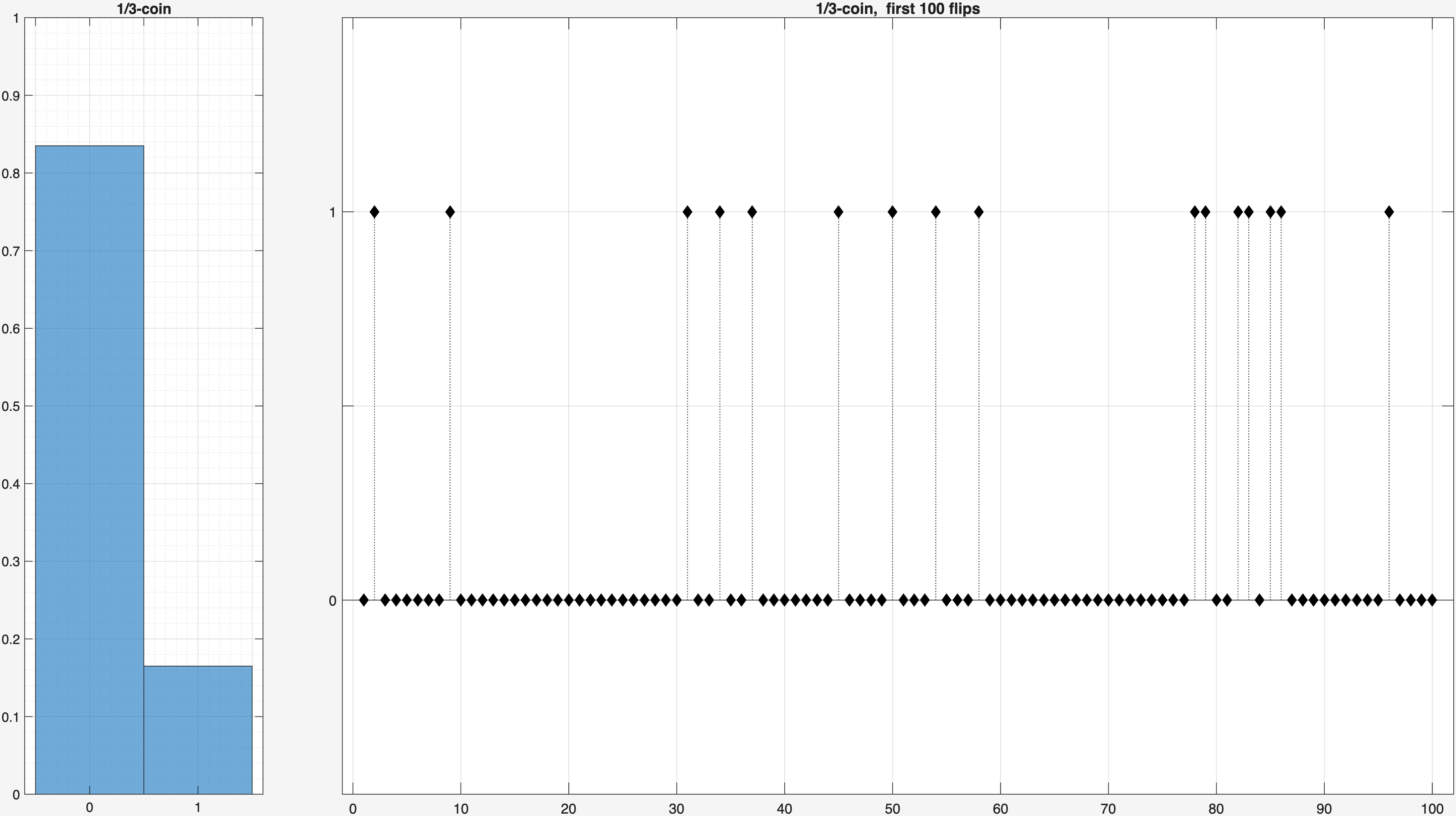}
\caption{A one-third-coin. \\
Flips: 10000; 
ones: 1650; 
zeros:  8350;  
expectation: 0.1650}
\end{figure}

\begin{figure}[p!]
\includegraphics[width=0.8\textwidth,height=0.24\textheight]{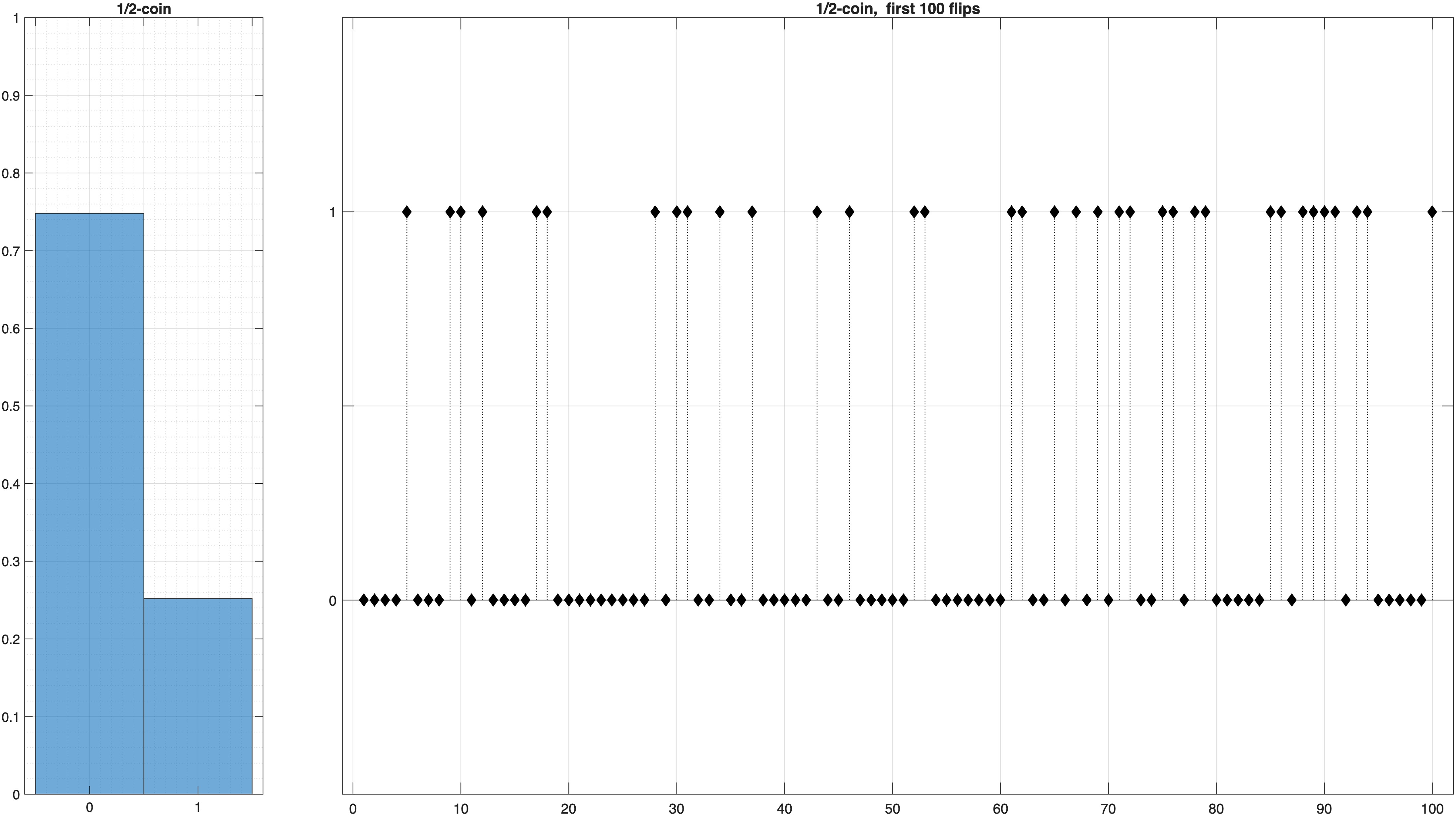}
\caption{A half-coin. \\
Flips: 10000; 
ones: 2519; 
zeros:  7481;  
expectation: 0.2519}
\end{figure}

\begin{figure}[p!]
\includegraphics[width=0.8\textwidth,height=0.24\textheight]{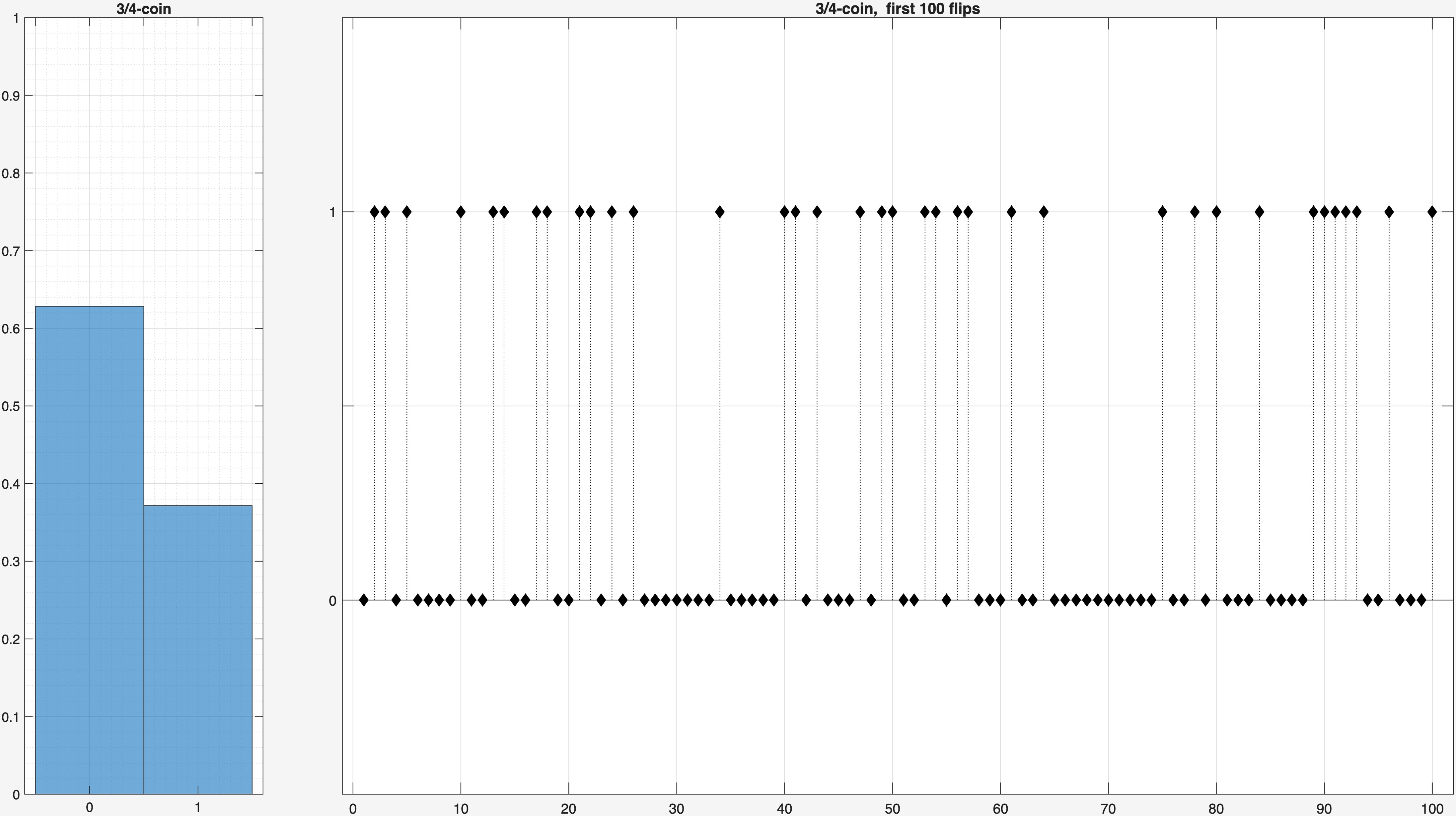}
\caption{A three-quarters-coin. \\
Flips: 10000; 
ones: 3716; 
zeros:  6284;  
expectation: 0.3716}
\end{figure}

\begin{figure}[p!]
\includegraphics[width=0.8\textwidth,height=0.24\textheight]{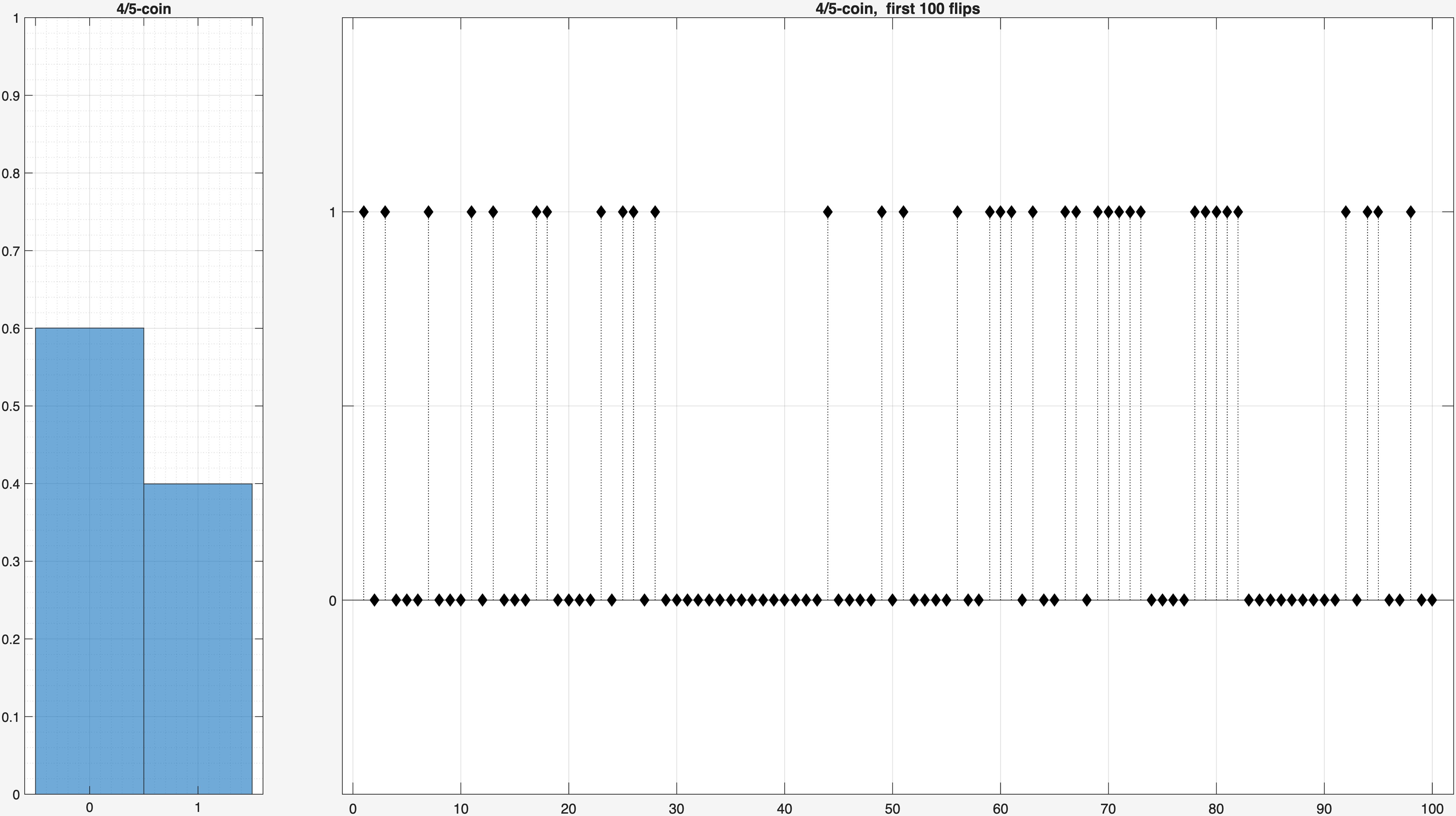}
\caption{A four-fifths-coin. \\
Flips: 10000; 
ones: 3997; 
zeros:  6003;  
expectation: 0.3997}
\end{figure}

\subsection{Examples of numerical simulation of tossing two partial coins}
\hspace*{\fill} \\
 
The function \texttt{twopartialcoins} of the dedicated MATLAB toolbox~\cite{Podlubny-partial-coin-toolbox} serves for the numerical simulation of tossing a pair of partial coins. 
Each partial coin is independently simulated as described in Section~\ref{sec:simulation},
and each flip of a pair of two partial coins produces the output equal to 0, 1, or 2.

\begin{figure}
\includegraphics[width=0.8\textwidth,height=0.24\textheight]{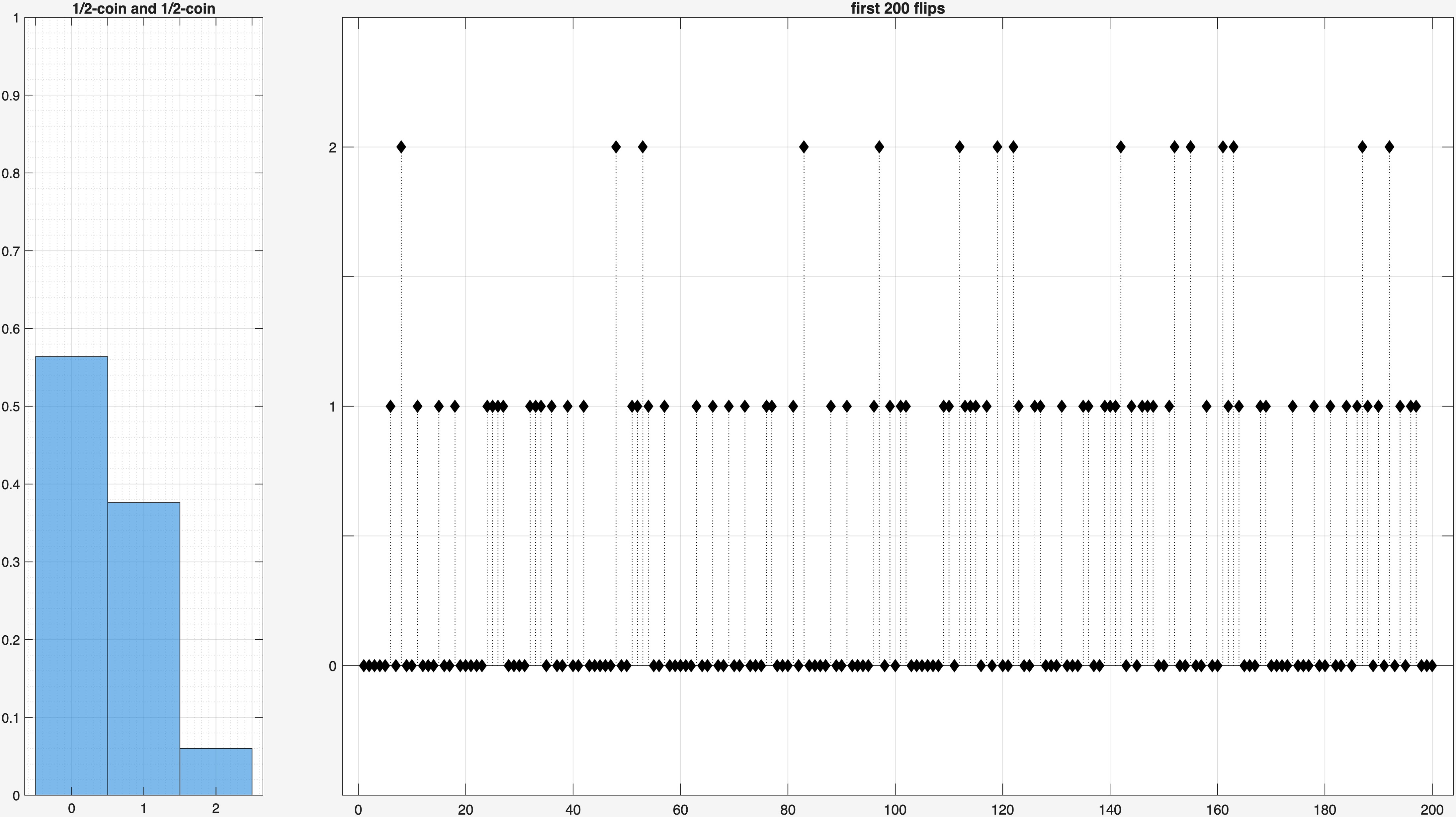}
\caption{Two half-coins. \\
Flips: 20000; 
ones on both partial coins: 1200;
ones on one partial coin: 7525; 
zeros on both partial coins:  11275;  
expectation: 0.4962} 
\label{fig:half-and-half-coins}
\end{figure}

\begin{figure}
\includegraphics[width=0.8\textwidth,height=0.24\textheight]{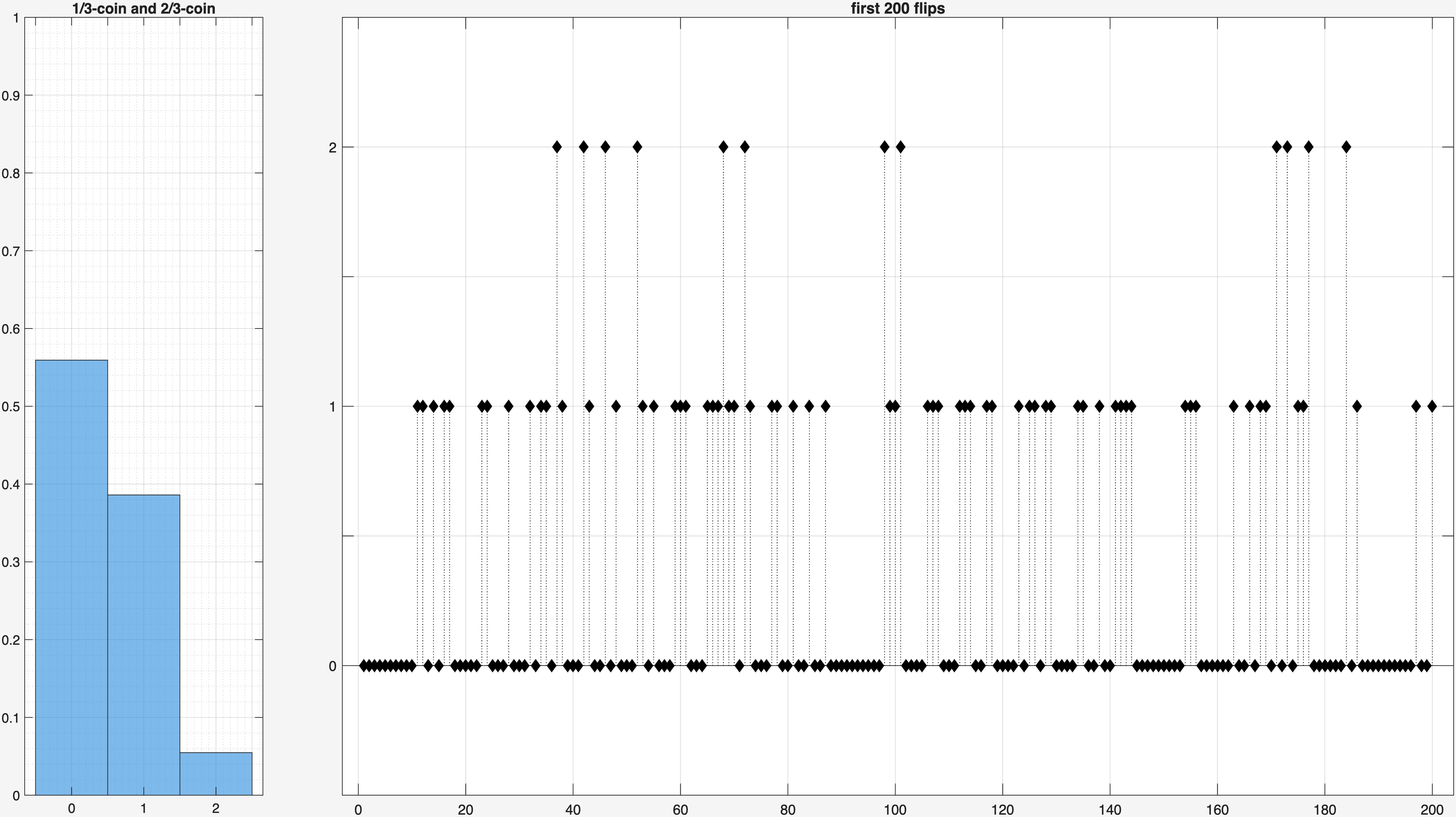}
\caption{A one-third-coin and a two-thirds-coin. \\
Flips: 20000; 
ones on both partial coins: 1029;
ones on one partial coin: 7793; 
zeros on both partial coins:  11178;  
expectation: 0.4925} 
\label{fig:one-third-and-two-thirds-coins}
\end{figure}

\begin{figure}
\includegraphics[width=0.8\textwidth,height=0.24\textheight]{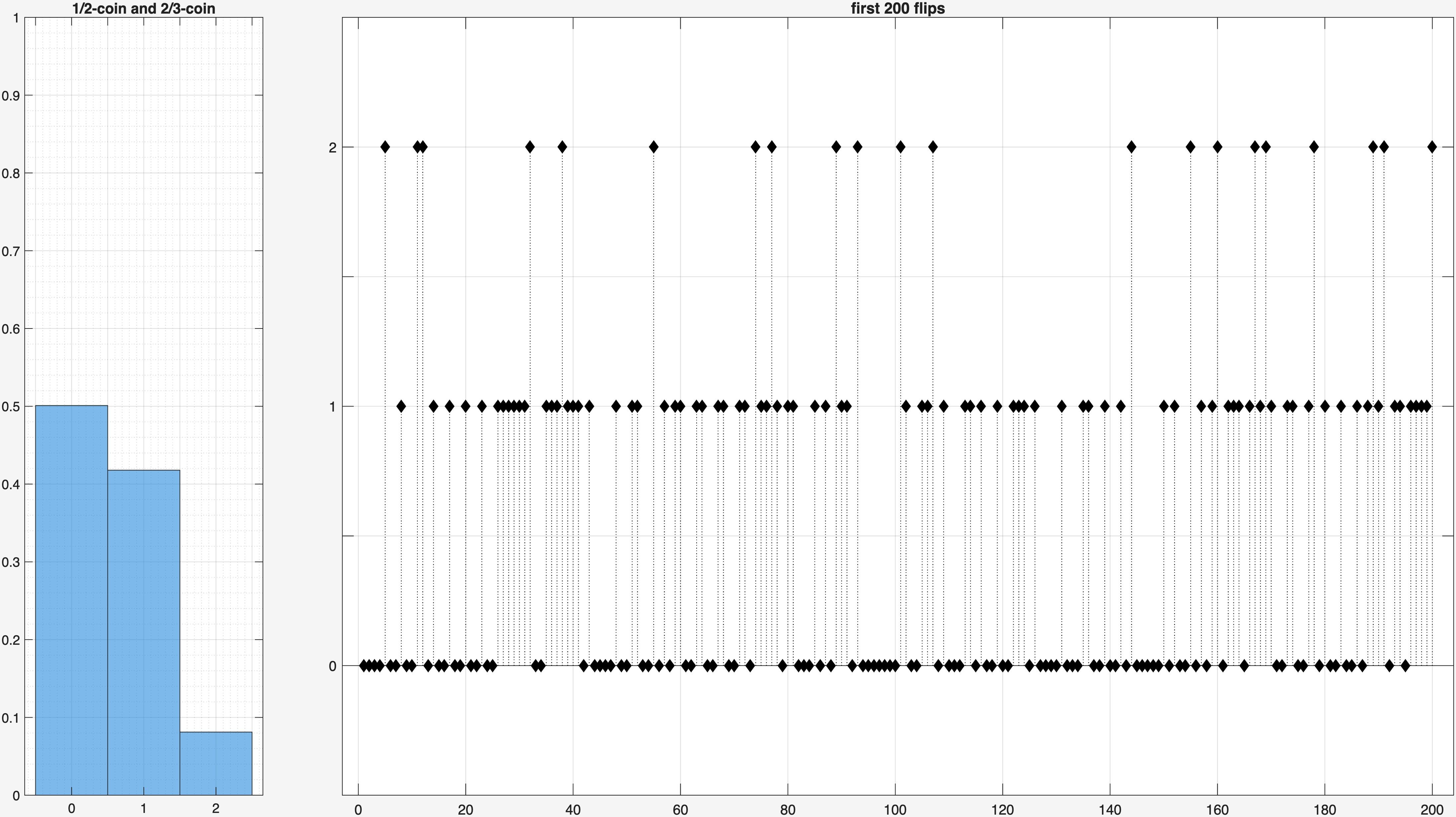}
\caption{A half-coin and a two-thirds-coin. \\
Flips: 20000; 
ones on both partial coins: 1624;
ones on one partial coin: 8357; 
zeros on both partial coins:  10019;  
expectation: 0.5802} 
\label{fig:half-and-two-thirds-coins}
\end{figure}

The results of the first two hundred flips of two partial coins are shown in the right 
of Figures~\ref{fig:half-and-half-coins}--\ref{fig:half-and-two-thirds-coins}, and the histograms are shown on the left. 

In Figure~\ref{fig:half-and-half-coins}, the results of the numerical simulation of flipping two half-coins are shown, and the expectation close to 0.5 is the same as in the case of flipping one fair coin. 

Flipping a one-third-coin and a a two-thirds-coin also gives the expectation close to 0.5, as in the case of flipping one fair coin (Figure~\ref{fig:one-third-and-two-thirds-coins}). 

In Figure~\ref{fig:half-and-two-thirds-coins}, the simulation results are shown for a half-coin and a two-thirds-coin.

\subsection{Tossing biased partial coins}
\hspace*{\fill} \\

The normal biased coin has unequal probabilities $a$ and $b=1-a$ for its sides, and its probability generating function is
\begin{equation}\label{eq:f-hat}
	\hat{f}(x) = a + bx
\end{equation}

Let us take $b<a$.
For a biased partial coin,  the probability generating function is 
\begin{equation} \label{eq:f-hat-series}
	(a+bx)^\mu = 
	a^\mu 
	\sum_{n=0}^{\infty}
		\left(\frac{b}{a}\right)^n 
		p_n x^n,
	\qquad
	p_n = {\mu \choose n},
\end{equation}
and the coefficients of this series have alternating signs 
exactly as in (\ref{eq:alternating-signs}).

In this case, the corresponding pairs
$\hat{g}_k$ and $\hat{h}_k$ ($ k = -1, 0, 1, 2, \ldots$)
 of nonnegative probability distributions for applying Theorem \ref{theorem-1} are defined by the following probability generating functions:
\begin{equation} \label{eq:g-hat}
	\hat{g}_k(x) = (1 - a^{-\mu} (a-bx)^\mu) \, x^k
\end{equation}

\begin{equation}  \label{eq:h-hat}
	\hat{h}_k(x) = ( (a+bx)^\mu - a^{-\mu} ( a^2 - b^2x^2  )^\mu )  \, x^k.
\end{equation} 

The procedure of simulation is the same as in Section~\ref{sec:simulation}. 

In the case of $a<b$, the variables $a$ and $b$ in relationships 
(\ref{eq:f-hat})--(\ref{eq:h-hat}) should be mutually swapped, 
and the simulated values of $\hat{F}_i$ (see step 5 section~\ref{sec:simulation})
are then computed as $\hat{F}_i = 1 -\hat{H}_i + \hat{G}_i$. 

The function \texttt{biasedpartialcoin}, which is also included in the dedicated MATLAB toolbox~\cite{Podlubny-partial-coin-toolbox}, serves for the numerical simulation of tossing a biased partial coin. 

\begin{figure}[p!]
\includegraphics[width=0.8\textwidth,height=0.24\textheight]{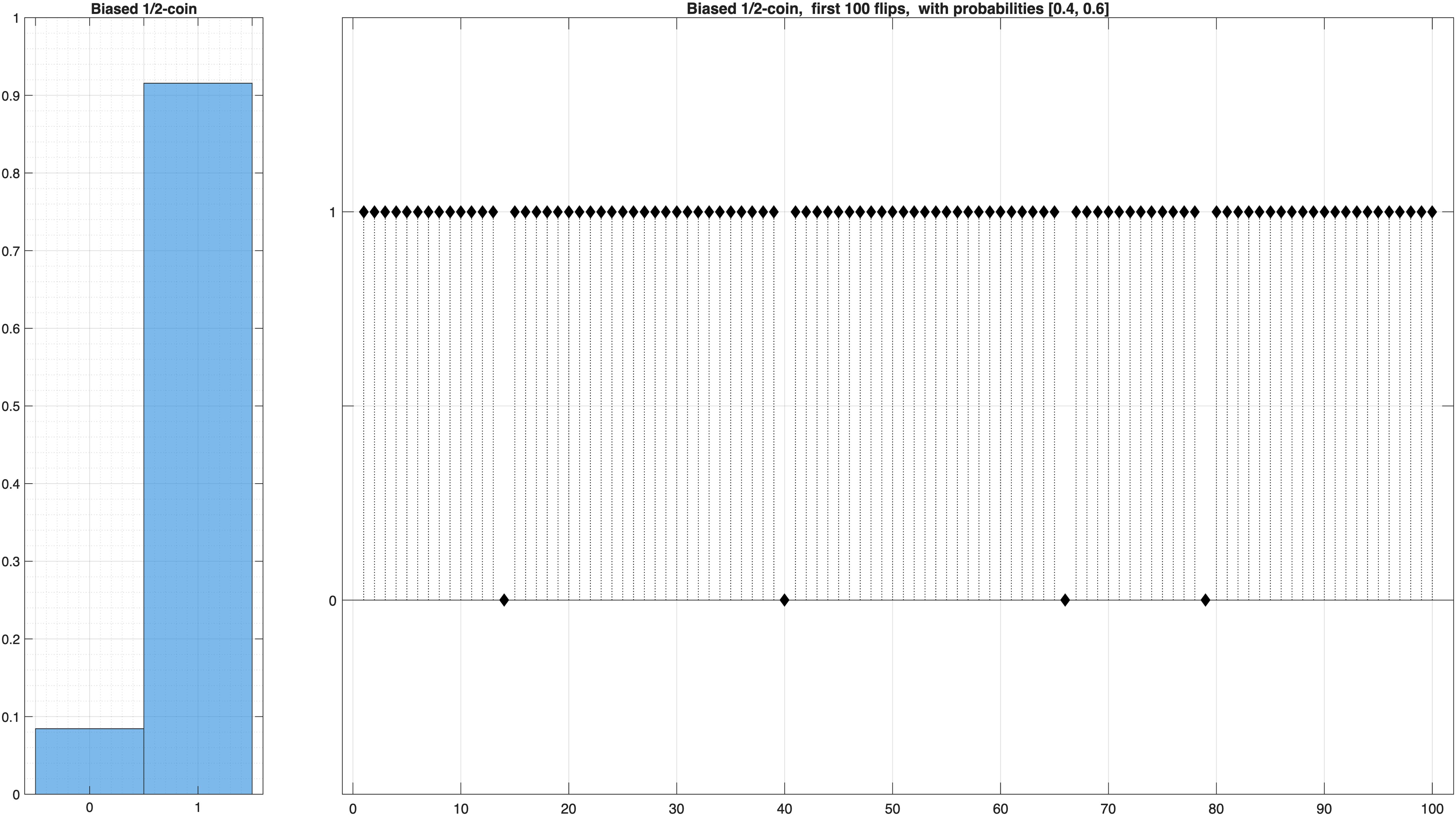}
\caption{A biased half-coin with $a = 0.4$ and $b=0.6$. \\
Flips: 10000; 
ones: 9156;
zeros:  844;  
expectation: 0.9156.}
\label{fig:biased-half-coin-04-06}
\end{figure}
\begin{figure}[p!]
\includegraphics[width=0.8\textwidth,height=0.24\textheight]{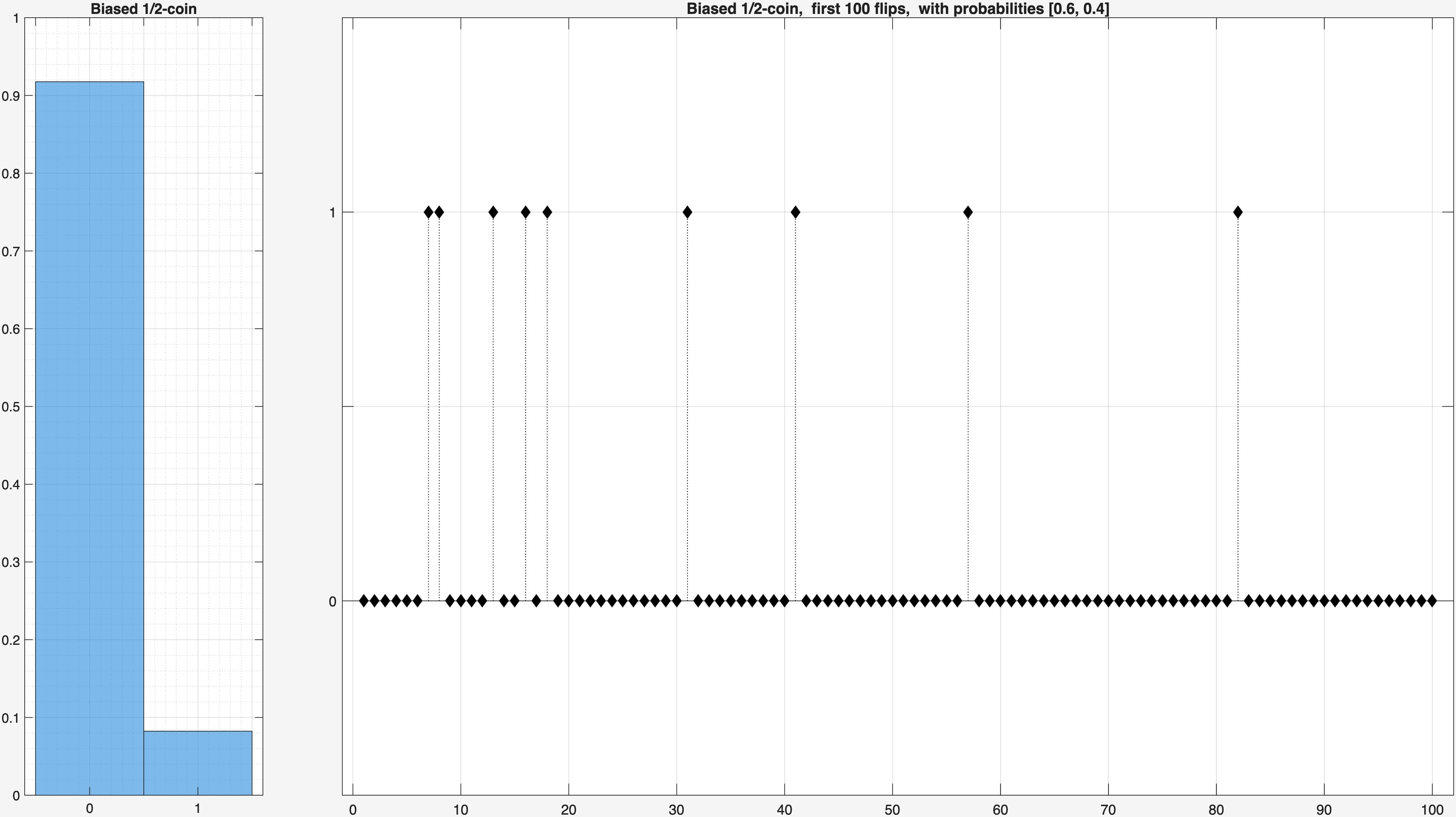}
\caption{A biased half-coin with $a = 0.6$ and $b=0.4$. \\
Flips: 10000; 
ones: 824;
zeros:  9176 ; 
expectation: 0.0824.} 
\label{fig:biased-half-coin-06-04}
\end{figure}
\begin{figure}[p!]
\includegraphics[width=0.8\textwidth,height=0.24\textheight]{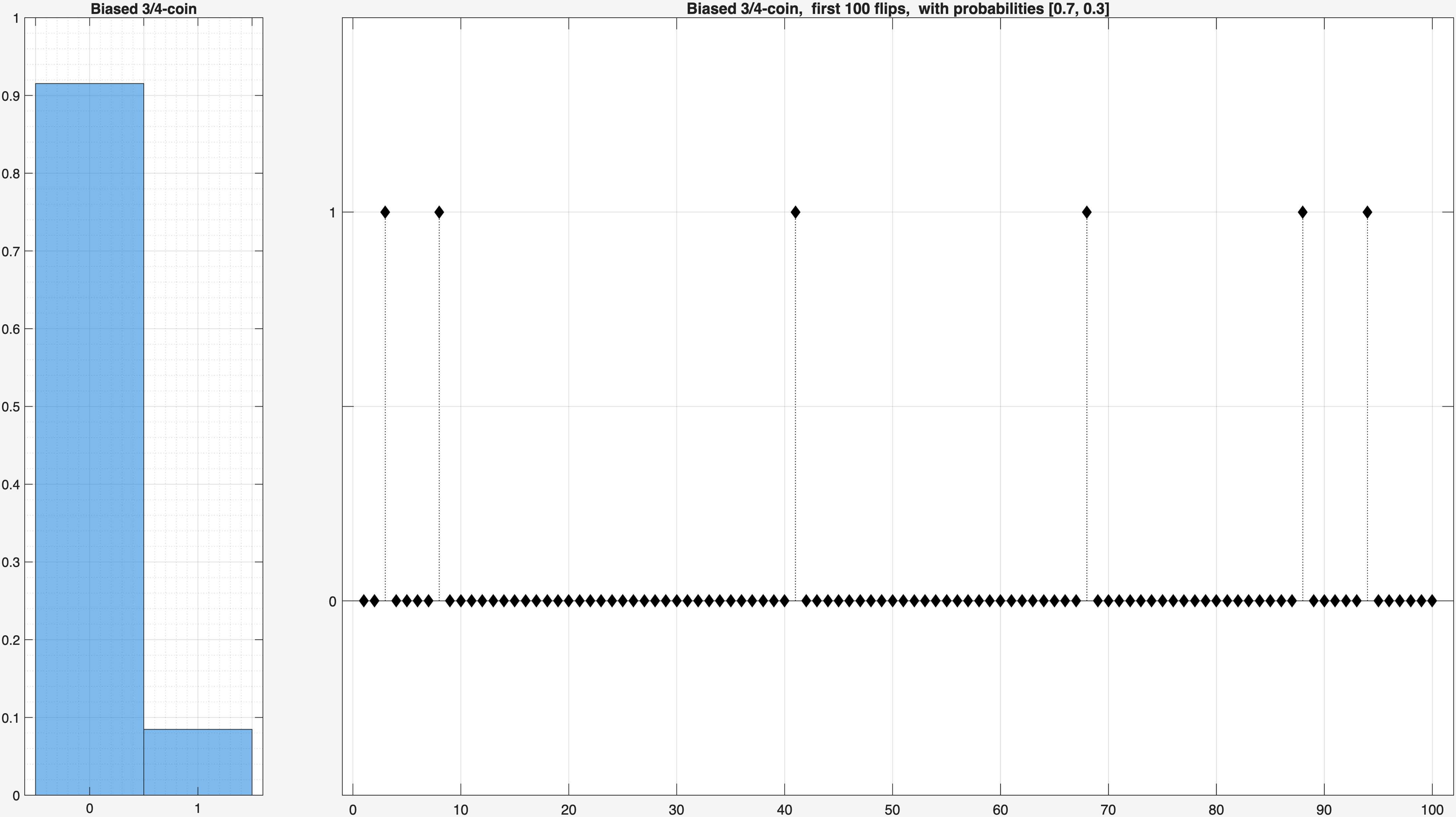}
\caption{A biased three-quarters-coin with $a = 0.7$ and $b=0.3$. \\
Flips: 10000 ; 
ones: 847;
zeros:  9153;  
expectation: 0.0847.} 
\label{fig:biased-3-4-coin-07-03}
\end{figure}

The results of simulations of a biased half-coin for $a=0.4, b=0.6$ and for $a=0.6, b=0.4$
are shown in Figure~\ref{fig:biased-half-coin-04-06} and in Figure~\ref{fig:biased-half-coin-06-04}, respectively. In Figure~\ref{fig:biased-3-4-coin-07-03} the results of simulation of a three-quarters-coin with $a = 0.7$ and $b=0.3$ are shown.

\section{Conclusion}

The case of a finite number of negative probabilities 
in a signed probability distribution
was considered and used for the Monte Carlo fractional differentiation 
in our recent works \cite{LP-MCFD-1,LP-MCFD-2,LP-MCFD-3}. 

To our best knowledge, the presented paper provides the first examples
of numerical simulation of signed probability distributions
with infinite number of negative probabilities.

The open question is how the probability generating functions $g(x)$ and $h(x)$
can be constructed in the case of a probability generating function $f(x)$ 
corresponding to a general case of a signed probability distribution. 

 \section{Acknowledgements}

\begin{figure}
\includegraphics[width=0.5\textwidth]{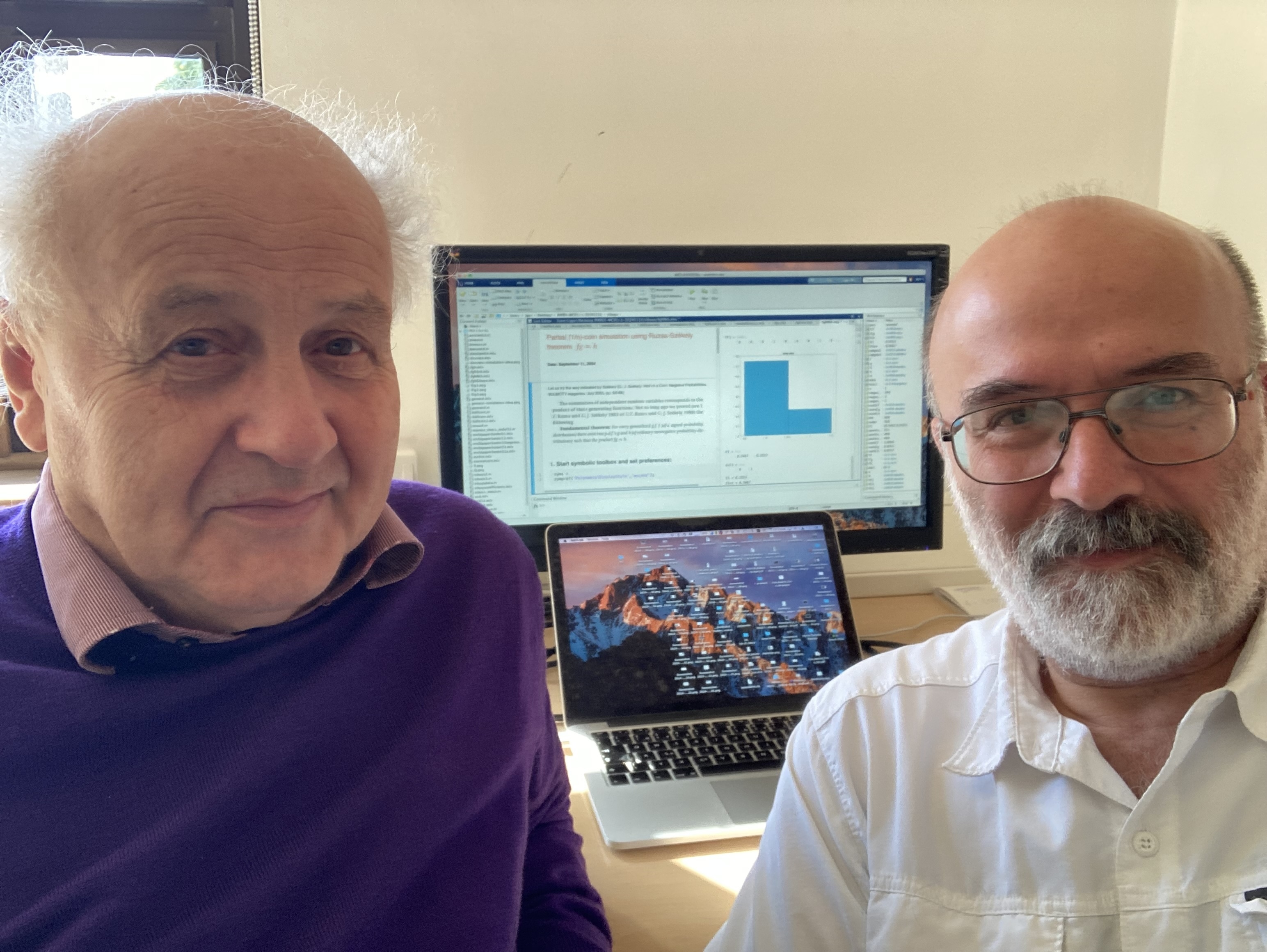}
\caption{September 13th, 2024, room M13, Isaac Newton Institute, Cambridge: 
the first successful simulation of a half-coin.} \label{fig:NL_IP_20240913}
\end{figure}

Nikolai Leonenko (NL) and Igor Podlubny (IP) would like to thank for support and hospitality during the programmes ``Fractional Differential Equations'' (FDE2),  ``Uncertainly Quantification and Modelling of Materials'' (USM), both supported by EPSRC grant EP/R014604/1, and the programme ``Stochastic systems for anomalous diffusion'' (SSD), supported by EPSRC grant EP/Z000580/1, at Isaac Newton Institute for Mathematical Sciences, Cambridge. The first successful simulation of a half-coin tossing happened during the SSD programme 
on September 13th, 2024, in the room M13 of the Isaac Newton Institute 
(Figure~\ref{fig:NL_IP_20240913}).

Also, NL was partially supported under the ARC Discovery Grant  DP220101680 (Australia), Croatian Scientific Foundation (HRZZ) grant ``Scaling in Stochastic Models'' (IP-2022-10-8081), grant FAPESP 22/09201-8 (Brazil) and the Taith Research Mobility grant (Wales, Cardiff University). 

The work of IP was also partially supported by grants VEGA 1/0674/23, APVV-18-0526, APVV-22-0508, and ARO W911NF-22-1-0264.

\end{document}